\documentclass[12pt]{article}
\usepackage[russian]{babel}
\usepackage[cp1251]{inputenc}
\usepackage{amsmath,latexsym}
\usepackage[psamsfonts]{amssymb}
\usepackage[dvips]{graphicx}
\usepackage[mathcal]{euscript}
\begin{document}
УДК 517.953
\begin{center}
\textbf{ ОБ УСЛОВИЯХ РАЗРЕШИМОСТИ КРАЕВЫХ ЗАДАЧ\\ ДЛЯ УРАВНЕНИЯ
ВЫСОКОГО ПОРЯДКА С РАЗРЫВНЫМ КОЭФФИЦИЕНТОМ}\\Б.Ю.Иргашев
  \\ Наманганский инженерно-строительный институт,\\
   Институт Математики им.Романовского АН.Р.Уз.
  \\ г. Наманган,Узбекистан
\\
  E-mail : bahromirgasev@gmail.com\end{center}

\textbf{Краткая анотация:} В работе расмотрена краевая задача для
уравнения Лаврентьева-Бицадзе высокого порядка.Найдены необходимые
и достаточные условия единственности решения.При обосновании
существования возникает проблема "малых знаменателей". Найдены
достаточные условия отделимости "малого знаменателя" от нуля.
Приведен пример краевой задачи не разрешимой методом Фурье, в случаи
прямоугольника с целочисленными измерениями сторон.

\textbf{Ключевые слова:} Дифференциальное уравнение, высокий
 порядок, уравнение Лаврентьева-Бицадзе, спектральный
 метод, единственность, существование, "малые
 знаменатели", бесконечный ряд, равномерная
 сходимость, алгебраическое число, неразрешимость.

\begin{center}
\textbf{1. Введение и формулировка основных результатов
}\end{center}

Рассмотрим уравнение в частных производных
\[Lu \equiv D_x^{2n}u\left( {x,y} \right) + \left( {{\mathop{\rm sgn}} y} \right)\,D_y^{2n}u\left( {x,y} \right) = 0,\eqno(1)\]
в прямоугольной области $\Omega  = \left\{ {\left( {x,y} \right):0 < x < l, - a < y < a} \right\}$, где $l,a$  - заданные положительные действительные числа,$n \in N$  . Пусть $\Omega _ + =
\Omega  \cap \left( {y > 0} \right),\,\,\Omega _ -   = \Omega \cap
\left( {y < 0} \right).$ Изучим для этого уравнения следующую краевую задачу.

\textbf{Задача D.} Найти в области $\Omega $  функцию $u(x,y)$  удовлетворяющую условиям:
\[
u \in C^{2n - 1} \left( {\overline \Omega  } \right) \cap C^{2n}
\left( {\Omega _ +   \cup \Omega _ -  } \right),\eqno(2)
\]
\[
Lu\left( {x,y} \right) \equiv 0,\,\,\,\,\left( {x,y} \right) \in
\Omega _ +   \cup \Omega _ -  ,\eqno(3)
\]
\[
D_x^{2s} u\left( {0,y} \right) = D_x^{2s} u\left( {l,y} \right) =
0,\,\,\, - a \le y \le a,\,\eqno(4)
\]
\[D_y^{q + \gamma s}u\left( {x, - a} \right) = {\varphi _s}\left( x \right),\,\,\,0 \le x \le l,\eqno(5)\]
\[D_y^{\chi  + \delta s}u\left( {x,a} \right) = {\psi _s}\left( x \right),\,\,0 \leqslant x \leqslant l,\eqno(6)
\]
\[\int\limits_0^l {{{\left[ {D_y^{2n}u\left( {x,y} \right)} \right]}^2}dx \leqslant M - const,\,\, - a \leqslant y \leqslant a,} \]
где ${\varphi _s}\left( x \right),{\psi _s}\left( x \right)$  -достаточно гладкие функции и для них выполняются условия согласования,$s = 0,...,n - 1;\,\,$  $\gamma ,\delta  \in \left\{ {1;2} \right\}$ если $\gamma  = \delta  = 1$, то $q,\chi  \in \left\{ {0,1,...,n} \right\};$ если $\gamma  = \delta  = 2$  , то   $q,\chi  \in \left\{ {0,1} \right\}.$

   При $n=1$  уравнение (1) есть известное уравнение Лаврентьева-Бицадзе, для которого некорректность задачи Дирихле было показано А.В.Бицадзе [1]. После этого специалистами, различными методами, были найдены условия единственности решения задачи Дирихле, как для уравнений 2-го порядка смешанного типа так и для уравнений высокого порядка с гладкими коэффициентами, например в работах [2]-[4]. Отметим, что в монографии Б.И.Пташника [4], имеется обширная литература по данной тематике.Из последних полученных результатов,
близких к изучаемой теме отметим работы [6]-[13].

В данной работе установлен критерий единственности решения задачи Дирихле. Решение построено в виде суммы по собственным функциям одномерной задачи. При обосновании  сходимости ряда возникает проблема малых знаменателей. Получены условия отделимости малого знаменателя от нуля.
\begin{center}
\textbf{2. Единственность решения
}\end{center}

Пусть  - решение уравнения (1) с условиями  (2)-(6). Рассмотрим его коэффициенты Фурье
\[{u_k}\left( y \right) = \int\limits_0^l {u\left( {x,y} \right){X_k}\left( x \right)dx} ,\,\,k = 1,2,...\,,\]
\[{X_k}\left( x \right) = \sqrt {\frac{2}{l}} \sin \frac{{\pi k}}{l}x.\eqno(7)\]
На основании (7) введем функции
\[{u_{k,\varepsilon }}\left( y \right) = \int\limits_\varepsilon ^{l - \varepsilon } {u\left( {x,y} \right){X_k}\left( x \right)dx} ,\eqno(8)\]
где $\epsilon >0$ - достаточно малое число. Дифференцируя равенство (8) по $y$  под знаком интеграла $2n$  раз, учитывая уравнение (1), при $y>0$  и $y<0$   , получим
\[u_{k,\varepsilon }^{\left( {2n} \right)}\left( y \right) = \int\limits_\varepsilon ^{l - \varepsilon } {D_y^{2n}u\left( {x,y} \right){X_k}\left( x \right)xdx}  =  - \int\limits_\varepsilon ^{l - \varepsilon } {D_x^{2n}u\left( {x,y} \right){X_k}\left( x \right)dx} ,\,\left( {y > 0} \right),\eqno(9)\]
\[u_{k,\varepsilon }^{\left( {2n} \right)}\left( y \right) = \int\limits_\varepsilon ^{l - \varepsilon } {D_y^{2n}u\left( {x,y} \right){X_k}\left( x \right)dx}  = \int\limits_\varepsilon ^{l - \varepsilon } {D_x^{2n}u\left( {x,y} \right){X_k}\left( x \right)dx} ,\,\left( {y < 0} \right),\eqno(10)\]
в интегралах из правых частей равенств (9) и (10), интегрируя по частям $2n$  раза и переходя к пределу при $\varepsilon  \to  + 0 $   с учетом граничных условий (4), получим обыкновенное  дифференциальное уравнение
\[u_k^{\left( {2n} \right)}\left( y \right) + {\left( { - 1} \right)^n}\left( {{\mathop{\rm sgn}} y} \right){\left( {\frac{{\pi k}}{l}} \right)^{2n}}{u_k}\left( y \right) = 0.\]
Рассмотрим отдельно случаи четного и нечетного $n$  . Пусть порядок уравнения равен $4m$  . Учитывая непрерывность функции $u(x,y)$  и её производных, на линии $y=0$ , получим следующую задачу:
\[\left\{ \begin{gathered}
  u_k^{\left( {4m} \right)}\left( y \right) + \operatorname{sgn} y{\left( {\frac{{\pi k}}{l}} \right)^{4m}}{u_k}\left( y \right) = 0, \hfill \\
  u_k^{\left( {q + \gamma j} \right)}\left( { - a} \right) = {\varphi _{jk}}, \hfill \\
  u_k^{\left( {\chi  + \delta s} \right)}\left( a \right) = {\psi _{jk}}, \hfill \\
  u_k^{\left( p \right)}\left( { + 0} \right) = u_k^{\left( p \right)}\left( { - 0} \right),\,j = \overline {0,2m - 1} ,\,\,p = \overline {0,(4m - 1)} . \hfill \\
\end{gathered}  \right.\eqno(11)\]
где
\[{\varphi _{jk}} = \int\limits_0^l {{\varphi _j}\left( x \right){X_k}\left( x \right)dx,\,\,} {\psi _{jk}} = \int\limits_0^l {{\psi _j}\left( x \right){X_k}\left( x \right)dx\,} .\]
При $y>0$  общее решение  уравнения (11) имеет вид
\[{u_k}\left( y \right) = \sum\limits_{p = 0}^{2m - 1} {\left( {c_p^1{Y_{1p}}(y) + c_p^2{Y_{2p}}\left( y \right)} \right)} ,\]
где
\[{Y_{1p}}\left( y \right) = {e^{{\alpha _p}y}}\cos {\beta _p}y,\,\,{Y_{2p}}\left( y \right) = {e^{{\alpha _p}y}}\sin {\beta _p}y,\]
\[{\left( {{Y_{1p}}\left( y \right)} \right)^{\left( q \right)}} = {\left( {\frac{{\pi k}}{l}} \right)^q}{e^{{\alpha _p}y}}\cos \left( {{\beta _p}y + q{\theta _p}} \right),\,\,\]
\[{\left( {{Y_{2p}}\left( y \right)} \right)^{\left( q \right)}} = {\left( {\frac{{\pi k}}{l}} \right)^q}{e^{{\alpha _p}y}}\sin \left( {{\beta _p}y + q{\theta _p}} \right),\]
 \[{\alpha _p} = \frac{{\pi k}}{l}\cos {\theta _p},{\beta _p} = \frac{{\pi k}}{l}\sin {\theta _p},\,{\theta _p} = \frac{\pi }{{4m}}\left( {1 + 2p} \right),\]
 \[p = \overline {0,(2m - 1)} ,\,q = \overline {0,4m - 1} ,\,\,{\alpha _p} > 0,\,p = \overline {0,(m - 1)} .\]
При $y<0$ имеем
\[{u_k}\left( y \right) = {d_0}{e^{\frac{{\pi k}}{l}y}} + \sum\limits_{s = 1}^{2m - 1} {{e^{{\mu _s}y}}\left( {d_s^1\cos {\nu _s}y + d_s^2\sin {\nu _s}y} \right)}  + {d_{2m}}{e^{ - \frac{{\pi k}}{l}y}},\]
\[u_k^{\left( {q + \gamma j} \right)}\left( y \right) = {\left( {\frac{{\pi k}}{l}} \right)^{\left( {q + \gamma j} \right)}}\left( {{d_0}{e^{\frac{{\pi ky}}{l}}} + } \right.\]
\[\left. { + \sum\limits_{s = 1}^{2m - 1} {{e^{{\mu _s}y}}\left( {d_s^1\cos \left( {{\nu _s}y + \left( {q + \gamma j} \right){\sigma _s}} \right) + d_s^2\sin \left( {{\nu _s}y + \left( {q + \gamma j} \right){\sigma _s}} \right)} \right)}  + {{\left( { - 1} \right)}^{\left( {q + \gamma j} \right)}}{d_{2m}}{e^{ - \frac{{\pi ky}}{l}}}} \right),\]
где
\[{\mu _s} = \frac{{\pi k}}{l}\cos {\sigma _s},{\nu _s} = \frac{{\pi k}}{l}\sin {\sigma _s},\,{\sigma _s} = \frac{{\pi s}}{{2m}},s = \overline {0,(2m - 1)} ,{\mu _s} > 0,\,s = \overline {0,(m - 1)} ,\,{\mu _m} = 0.\]
Удовлетворив краевым условиям задачи (11) получим систему алгебраических уравнений
\[\left\{ \begin{gathered}
  \sum\limits_{p = 0}^{2m - 1} {{e^{{\alpha _p}y}}\left( {c_p^1\cos \left( {{\beta _p}a + \left( {\chi  + \delta j} \right){\theta _p}} \right) + c_p^2sin\left( {{\beta _p}a + \left( {\chi  + \delta j} \right){\theta _p}} \right)} \right)}  =  \hfill \\
   = {\left( {\frac{l}{{\pi k}}} \right)^{\chi  + \delta j}}{\psi _{jk}}, \hfill \\
  \sum\limits_{s = 0}^{2m - 1} {{e^{ - a{\mu _s}}}\left( {d_s^1\cos \left( { - {\nu _s}a + \left( {q + \gamma j} \right){\sigma _s})} \right)} \right. + }  \hfill \\
   + \left. {d_s^2\sin \left( { - {\nu _s}a + \left( {q + \gamma j} \right){\sigma _s}} \right)} \right) = {\left( {\frac{l}{{\pi k}}} \right)^{q + j\gamma }}{\varphi _{jk}}, \hfill \\
  \sum\limits_{p = 0}^{2m - 1} {\left( {c_p^1\cos \left( {t{\theta _p}} \right) + c_p^2\sin \left( {t{\theta _p}} \right)} \right)}  =  \hfill \\
   = {d_0} + \sum\limits_{p = 1}^{2m - 1} {\left( {d_p^1\cos \left( {t{\sigma _p}} \right) + d_p^2\sin \left( {t{\sigma _p}} \right)} \right) + {{\left( { - 1} \right)}^t}{d_{2m}},}  \hfill \\
  \,j = \overline {0,2m - 1} ,\,\,\,t = \overline {0,(4m - 1)} . \hfill \\
\end{gathered}  \right.\eqno(12)\]
Введем следующие обозначения :
\[{\omega _{j,p}} = {\beta _p}a + \left( {\chi  + \delta j} \right){\theta _p},\,{\tau _{j,s}} =  - {\nu _s}a + \left( {q + \gamma j} \right){\sigma _s},\]
\[A_{2m,2m}^ +  = \left( {\begin{array}{*{20}{c}}
  {{e^{{\alpha _0}a}}\cos {\omega _{0,0}}}&{{e^{{\alpha _0}a}}\sin {\omega _{0,0}}}&.&{{e^{{\alpha _{m - 1}}a}}\cos {\omega _{0,m - 1}}}&{{e^{{\alpha _{m - 1}}a}}\sin {\omega _{0,m - 1}}} \\
  {{e^{{\alpha _0}a}}\cos {\omega _{1,0}}}&{{e^{{\alpha _0}a}}\sin {\omega _{1,0}}}&.&{{e^{{\alpha _{m - 1}}a}}\cos {\omega _{1,m - 1}}}&{{e^{{\alpha _{m - 1}}a}}\sin {\omega _{1,m - 1}}} \\
  .&.&.&.&. \\
  {{e^{{\alpha _0}a}}\cos {\omega _{2m - 1,0}}}&{{e^{{\alpha _0}a}}\sin {\omega _{2m - 1,0}}}&.&{{e^{{\alpha _0}a}}\cos {\omega _{2m - 1,m - 1}}}&{{e^{{\alpha _{m - 1}}a}}\sin {\omega _{2m - 1,m - 1}}}
\end{array}} \right),\]
\[A_{2m,2m}^ -  = \left( {\begin{array}{*{20}{c}}
  {{e^{{\alpha _m}a}}\cos {\omega _{0,m}}}&{{e^{{\alpha _m}a}}\sin {\omega _{0,m}}}&.&{{e^{{\alpha _{m2 - 1}}a}}\cos {\omega _{0,2m - 1}}}&{{e^{{\alpha _{2m - 1}}a}}\sin {\omega _{0,2m - 1}}} \\
  {{e^{{\alpha _m}a}}\cos {\omega _{1,m}}}&{{e^{{\alpha _m}a}}\sin {\omega _{1,m}}}&.&{{e^{{\alpha _{2m - 1}}a}}\cos {\omega _{1,2m - 1}}}&{{e^{{\alpha _{2m - 1}}a}}\sin {\omega _{1,2m - 1}}} \\
  .&.&.&.&. \\
  {{e^{{\alpha _m}a}}\cos {\omega _{2m - 1,m}}}&{{e^{{\alpha _m}a}}\sin {\omega _{2m - 1,m}}}&.&{{e^{{\alpha _{2m - 1}}a}}\cos {\omega _{2m - 1,2m - 1}}}&{{e^{{\alpha _{2m - 1}}a}}\sin {\omega _{2m - 1,2m - 1}}}
\end{array}} \right),\]
\[B_{2m,2m - 1}^ +  = \left( {\begin{array}{*{20}{c}}
  {{e^{ - {\mu _{m + 1}}a}}\cos {\tau _{0,m + 1}}}&.&.&{{e^{ - {\mu _{2m - 1}}a}}\sin {\tau _{0,2m - 1}}}&{{{\left( { - 1} \right)}^q}{e^{\frac{{\pi k}}{l}a}}} \\
  {{e^{ - {\mu _{m + 1}}a}}\cos {\tau _{1, + 1}}}&.&.&{{e^{ - {\mu _{2m - 1}}a}}\sin {\tau _{1,2m - 1}}}&{{{\left( { - 1} \right)}^{q + \gamma }}{e^{\frac{{\pi k}}{l}a}}} \\
  .&.&.&.&. \\
  {{e^{ - {\mu _{m + 1}}a}}\cos {\tau _{2m - 1,m + 1}}}&.&.&{{e^{ - {\mu _{2m - 1}}a}}\sin {\tau _{2m - 1,2m - 1}}}&{{{\left( { - 1} \right)}^{q + \gamma \left( {2m - 1} \right)}}{e^{\frac{{\pi k}}{l}a}}}
\end{array}} \right),\]
или  в компактной записи
\[B_{2m,2m - 1}^ +  = \left( {{e^{ - {\mu _s}a}}\cos {\tau _{j,s}},{e^{ - {\mu _s}a}}\sin {\tau _{j,s}},{{\left( { - 1} \right)}^{q + \gamma j}}{e^{\frac{{\pi k}}{l}a}}} \right)_{s = \overline {m + 1,2m - 1} }^{j = \overline {0,2m - 1} },\]
аналогично
\[B_{2m,2m - 1}^ -  = \left( {{e^{ - \frac{{\pi k}}{l}a}},{e^{ - {\mu _s}a}}\cos {\tau _{j,s}},{e^{ - {\mu _s}a}}\sin {\tau _{j,s}}} \right)_{s = \overline {1,m - 1} }^{j = \overline {0,2m - 1} },\]
\[C_{4m,2m}^ +  = \left( {\cos j{\theta _s},\sin j{\theta _s}} \right)_{s = \overline {0,m - 1} }^{j = \overline {0,4m - 1} },\,C_{4m,2m}^ -  = \left( {\cos j{\theta _s},\sin j{\theta _s}} \right)_{s = \overline {m,2m - 1} }^{j = \overline {0,4m - 1} },\]
\[D_{4m,2m - 1}^ +  = \left( { - \cos j{\sigma _s}, - \sin j{\sigma _s},{{\left( { - 1} \right)}^j}} \right)_{s = \overline {m + 1,2m - 1} }^{j = \overline {0,4m - 1} },\]
\[D_{4m,2m - 1}^ -  = \left( {1, - \cos j{\sigma _s}, - \sin j{\sigma _s}} \right)_{s = \overline {1,m - 1} }^{j = \overline {0,4m - 1} }.\]
Используя формулу Эйлера ${e^{iz}} = \cos z + i\sin z$, сделаем некоторые формальные преобразования
\[B_{2m,2}^0 = \left( {\begin{array}{*{20}{c}}
  {\cos {\tau _{0,m}}}&{\sin {\tau _{0,m}}} \\
  {\cos {\tau _{1,m}}}&{\sin {\tau _{1,m}}} \\
  .&. \\
  {\cos {\tau _{2m - 1,m}}}&{\sin {\tau _{2m - 1,m}}}
\end{array}} \right) = \frac{i}{2}\left( {\begin{array}{*{20}{c}}
  {{e^{i{\tau _{0,m}}}}}&{{e^{ - i{\tau _{0,m}}}}} \\
  .&. \\
  {{e^{i{\tau _{2m - 1,m}}}}}&{{e^{ - i{\tau _{2m - 1,m}}}}}
\end{array}} \right) = \frac{i}{2}\left( {\begin{array}{*{20}{c}}
  {B_{2m,1}^ + }&{B_{2m,1}^ - }
\end{array}} \right),\]
\[D_{4m,2}^0 = \left( {\begin{array}{*{20}{c}}
  {\cos 0}&{\sin 0} \\
  {\cos \frac{\pi }{2}}&{\sin \frac{\pi }{2}} \\
  .&. \\
  {\cos \left( {4m - 1} \right)\frac{\pi }{2}}&{\sin \left( {4m - 1} \right)\frac{\pi }{2}}
\end{array}} \right) = \]
\[ = \frac{i}{2}\left( {\begin{array}{*{20}{c}}
  {{e^{i0 \cdot \frac{\pi }{2}}}}&{{e^{ - i0 \cdot \frac{\pi }{2}}}} \\
  .&. \\
  {{e^{i\left( {4m - 1} \right)\frac{\pi }{2}}}}&{{e^{ - i\left( {4m - 1} \right)\frac{\pi }{2}}}}
\end{array}} \right) = \frac{i}{2}\left( {\begin{array}{*{20}{c}}
  {D_{4m,1}^ + }&{D_{4m,1}^ - }
\end{array}} \right),\]
тогда основной определитель, системы (12) , будет иметь вид
\[{\Delta _{1k}} = \frac{i}{2}\det \left( {\begin{array}{*{20}{c}}
  {A_{2m,2m}^ + }&{A_{2m,2m}^ - }&0&0&0&0 \\
  0&0&{B_{2m,2m - 1}^ - }&{B_{2m,2m - 1}^ + }&{B_{2m,1}^ + }&{B_{2m,1}^ - } \\
  {C_{4m,2m}^ + }&{C_{4m,2m}^ - }&{D_{4m,2m - 1}^ - }&{D_{4m,2m - 1}^ + }&{D_{4m,1}^ + }&{D_{4m,1}^ - }
\end{array}} \right).\]
Найдем асимптотику определителя ${\Delta _{1k}}$ , при больших значениях $k$. Для этого вычислим слагаемое, куда входит экспонента с наибольшей положительной степенью. С точностью до знака он имеет вид
\[{\Delta _{2k}} = \frac{i}{2}\left| {A_{2m,2m}^ + } \right|\left( {\left| {B_{2m,2m - 1}^ + B_{2m,1}^ + } \right| \cdot \left| {C_{4m,2m}^ - D_{4m,2m - 1}^ - D_{4m,1}^ - } \right| - } \right.\]
\[ - \left. {\left| {B_{2m,2m - 1}^ + B_{2m,1}^ - } \right| \cdot \left| {C_{4m,2m}^ - D_{4m,2m - 1}^ - D_{4m,1}^ + } \right|} \right).\]
Перейдем к вычислениям. Имеем
\[\left| {A_{2m,2m}^ + } \right| = {\left( {\frac{i}{2}} \right)^m}{e^{2\alpha a}}\left| {\begin{array}{*{20}{c}}
  {{e^{i{\omega _{0,0}}}}}&{{e^{ - i{\omega _{0,0}}}}}&.&{{e^{i{\omega _{0,m - 1}}}}}&{{e^{ - i{\omega _{0,m - 1}}}}} \\
  {{e^{i{\omega _{1,0}}}}}&{{e^{ - i{\omega _{1,0}}}}}&.&{{e^{i{\omega _{1,m - 1}}}}}&{{e^{ - i{\omega _{1,m - 1}}}}} \\
  .&.&.&.&. \\
  {{e^{i{\omega _{2m - 1,0}}}}}&{{e^{ - i{\omega _{2m - 1,0}}}}}&.&{{e^{i{\omega _{2m - 1,m - 1}}}}}&{{e^{ - i{\omega _{2m - 1,m - 1}}}}}
\end{array}} \right| = \]
\[ = {\left( {\frac{i}{2}} \right)^m}{e^{2\alpha a}} \cdot \]
\[\left| {\begin{array}{*{20}{c}}
  {{e^{i\left( {{\beta _0}a + \chi {\theta _0}} \right)}}}&{{e^{ - i\left( {{\beta _0}a + \chi {\theta _0}} \right)}}}&.&{{e^{i\left( {{\beta _{m - 1}}a + \chi {\theta _{m - 1}}} \right)}}}&{{e^{ - i\left( {{\beta _{m - 1}}a + \chi {\theta _{m - 1}}} \right)}}} \\
  {{e^{i\left( {{\beta _0}a + \left( {\chi  + \delta } \right){\theta _0}} \right)}}}&{{e^{ - i\left( {{\beta _0}a + \left( {\chi  + \delta } \right){\theta _0}} \right)}}}&.&{{e^{i\left( {{\beta _{m - 1}}a + \left( {\chi  + \delta } \right){\theta _{m - 1}}} \right)}}}&{{e^{ - i\left( {{\beta _{m - 1}}a + \left( {\chi  + \delta } \right){\theta _{m - 1}}} \right)}}} \\
  .&.&.&.&. \\
  {{e^{i\left( {{\beta _0}a + \left( {\chi  + \delta \left( {2m - 1} \right)} \right){\theta _0}} \right)}}}&{{e^{ - i\left( {{\beta _0}a + \left( {\chi  + \delta \left( {2m - 1} \right)} \right){\theta _0}} \right)}}}&.&{{e^{i\left( {{\beta _{m - 1}}a + \left( {\chi  + \delta \left( {2m - 1} \right)} \right){\theta _{m - 1}}} \right)}}}&{{e^{ - i\left( {{\beta _{m - 1}}a + \left( {\chi  + \delta \left( {2m - 1} \right)} \right){\theta _{m - 1}}} \right)}}}
\end{array}} \right| = \]
\[ = {\left( {\frac{i}{2}} \right)^m}{e^{2\alpha a}}\left| {\begin{array}{*{20}{c}}
  1&1&.&1&1 \\
  {{e^{i\delta {\theta _0}}}}&{{e^{ - i\delta {\theta _0}}}}&.&{{e^{i\delta {\theta _{m - 1}}}}}&{{e^{ - i\delta {\theta _{m - 1}}}}} \\
  .&.&.&.&. \\
  {{e^{i\delta \left( {2m - 1} \right){\theta _0}}}}&{{e^{ - i\delta \left( {2m - 1} \right){\theta _0}}}}&.&{{e^{i\delta \left( {2m - 1} \right){\theta _{m - 1}}}}}&{{e^{ - i\delta \left( {2m - 1} \right){\theta _{m - 1}}}}}
\end{array}} \right| = \]
\[ = {\left( {\frac{i}{2}} \right)^m}{e^{2\alpha a}}\prod\limits_{j = 0}^{m - 1} {\left( { - 2i\sin \delta {\theta _j}} \right)} \prod\limits_{0 = s < j = m - 1}^{} {4\left( {1 - \cos \delta \left( {{\theta _j} - {\theta _s}} \right)} \right)\left( {1 - \cos \delta \left( {{\theta _j} + {\theta _s}} \right)} \right)}  \ne 0,\]
т.к.
\[0 < {\theta _j} < \frac{\pi }{2},\,\forall j,\]
\[{\theta _j} + {\theta _s} = \frac{\pi }{{4m}}\left( {1 + 2j} \right) + \frac{\pi }{{4m}}\left( {1 + 2s} \right) = \frac{{1 + j + s}}{{2m}}\pi  \leqslant \frac{{2m - 1}}{{2m}}\pi ,\]
здесь
\[\alpha  = {\alpha _0} + {\alpha _1} + ... + {\alpha _{m - 1}}.\]
Далее вводя обозначение $\mu  = {\mu _{m + 1}} + ... + {\mu _{2m - 1}}$  , имеем
\[\det \left( {B_{2m,2m - 1}^ + ,B_{2m,1}^ + } \right) = \det \left( {{e^{ - {\mu _s}a}}\cos {\tau _{j,s}},{e^{ - {\mu _s}a}}\sin {\tau _{j,s}},{{\left( { - 1} \right)}^{q + \gamma j}}{e^{\frac{{\pi k}}{l}a}},{e^{i{\tau _{j,m}}}}} \right)_{s = \overline {m + 1,2m - 1} }^{j = \overline {0,2m - 1} } = \]
\[ = {e^{ - 2a\mu  + \frac{{\pi k}}{l}a}}\det \left( {\cos {\tau _{j,s}},\sin {\tau _{j,s}},{{\left( { - 1} \right)}^{q + \gamma j}},{e^{i\left( { - {\nu _m}a + \left( {q + \gamma j} \right){\sigma _m}} \right)}}} \right)_{s = \overline {m + 1,2m - 1} }^{j = \overline {0,2m - 1} } = \]
\[ = {e^{ - 2a\mu  + \frac{{\pi k}}{l}a}}{\left( {\frac{i}{2}} \right)^{m - 1}}\det \left( {{e^{i{\tau _{j,s}}}},{e^{ - i{\tau _{j,s}}}},{{\left( { - 1} \right)}^{q + \gamma j}},{e^{i\left( { - \frac{{\pi ka}}{l} + \left( {q + \gamma j} \right)\frac{\pi }{2}} \right)}}} \right)_{s = \overline {m + 1,2m - 1} }^{j = \overline {0,2m - 1} } = \]
\[ = {e^{ - 2a\mu  + \frac{{\pi k}}{l}a}}{e^{ - \frac{{\pi ka}}{l}i}}{\left( {\frac{i}{2}} \right)^{m - 1}}\det \left( {{e^{i\left( { - {\nu _s}a + q{\sigma _s} + \gamma j{\sigma _s}} \right)}},{e^{ - i\left( { - {\nu _s}a + q{\sigma _s} + \gamma j{\sigma _s}} \right)}},{{\left( { - 1} \right)}^{q + \gamma j}},{i^{q + \gamma j}}} \right)_{s = \overline {m + 1,2m - 1} }^{j = \overline {0,2m - 1} } = \]
\[ = {e^{ - 2a\mu  + \frac{{\pi k}}{l}a}}{e^{ - \frac{{\pi ka}}{l}i}}{\left( {\frac{i}{2}} \right)^{m - 1}}{\left( { - i} \right)^q}\det \left( {{e^{i\gamma j{\sigma _s}}},{e^{ - i\gamma j{\sigma _s}}},{{\left( { - 1} \right)}^{\gamma j}},{i^{\gamma j}}} \right)_{s = \overline {m + 1,2m - 1} }^{j = \overline {0,2m - 1} } = \]
\[ = {e^{ - 2a\mu  + \frac{{\pi k}}{l}a}}{e^{ - \frac{{\pi ka}}{l}i}}{\left( {\frac{i}{2}} \right)^{m - 1}}{\left( { - i} \right)^q}{M_1}\left( {{i^\gamma } - {{\left( { - 1} \right)}^\gamma }} \right)\prod\limits_{s = m + 1}^{2m - 1} {\left( {{i^\gamma } - {e^{i\gamma {\sigma _s}}}} \right)\left( {{i^\gamma } - {e^{ - i\gamma {\sigma _s}}}} \right)}  = \]
\[ = {e^{ - 2a\mu  + \frac{{\pi k}}{l}a}}{e^{ - \frac{{\pi ka}}{l}i}}{\left( {\frac{i}{2}} \right)^{m - 1}}{\left( { - i} \right)^q}{M_1}\left( {{i^\gamma } - {{\left( { - 1} \right)}^\gamma }} \right)\prod\limits_{s = m + 1}^{2m - 1} {\left( {{{\left( { - 1} \right)}^\gamma } + 1 - 2{i^\gamma }\cos \gamma {\sigma _s}} \right)} ,\]
где
\[{M_1} = \det \left( {{e^{i\gamma j{\sigma _s}}},{e^{ - i\gamma j{\sigma _s}}},{{\left( { - 1} \right)}^{\gamma j}}} \right)_{s = \overline {m + 1,2m - 1} }^{j = \overline {0,2m - 1} } \ne 0,\]
далее
\[\det \left( {B_{2m,2m - 1}^ + ,B_{2m,1}^ - } \right) = \det \left( {{e^{ - {\mu _s}a}}\cos {\tau _{j,s}},{e^{ - {\mu _s}a}}\sin {\tau _{j,s}},{{\left( { - 1} \right)}^{q + \gamma j}}{e^{\frac{{\pi k}}{l}a}},{e^{ - i{\tau _{j,m}}}}} \right)_{s = \overline {m + 1,2m - 1} }^{j = \overline {0,2m - 1} } = \]
\[ = {e^{ - 2\mu a + \frac{{\pi k}}{l}a}}\det \left( {\cos {\tau _{j,s}},\sin {\tau _{j,s}},{{\left( { - 1} \right)}^{q + \gamma j}},{e^{ - i{\tau _{j,m}}}}} \right)_{s = \overline {m + 1,2m - 1} }^{j = \overline {0,2m - 1} } = \]
\[ = {e^{ - 2\mu a + \frac{{\pi k}}{l}a}}{\left( {\frac{i}{2}} \right)^{m - 1}}{\left( { - 1} \right)^q} \cdot \]
\[ \cdot \det \left( {{e^{i\left( { - {\nu _s}a + \left( {q + \gamma j} \right){\sigma _s}} \right)}},{e^{ - i\left( { - {\nu _s}a + \left( {q + \gamma j} \right){\sigma _s}} \right)}},{{\left( { - 1} \right)}^{\gamma j}},{e^{ - i\left( { - \frac{{\pi k}}{l}a + \left( {q + \gamma j} \right)\frac{\pi }{2}} \right)}}} \right)_{s = \overline {m + 1,2m - 1} }^{j = \overline {0,2m - 1} } = \]
\[ = {e^{ - 2\mu a + \frac{{\pi k}}{l}a}}{e^{i\frac{{\pi k}}{l}a}}{\left( {\frac{i}{2}} \right)^{m - 1}}{i^q}\det \left( {{e^{i\gamma j{\sigma _s}}},{e^{ - i\gamma j{\sigma _s}}},{{\left( { - 1} \right)}^{\gamma j}},{{\left( { - i} \right)}^{j\gamma }}} \right)_{s = \overline {m + 1,2m - 1} }^{j = \overline {0,2m - 1} } = \]
\[ = {e^{ - 2\mu a + \frac{{\pi k}}{l}a}}{e^{i\frac{{\pi k}}{l}a}}{\left( {\frac{i}{2}} \right)^{m - 1}}{i^q}\left( {{{\left( { - i} \right)}^\gamma } - {{\left( { - 1} \right)}^\gamma }} \right){M_1}\prod\limits_{s = m + 1}^{2m - 1} {\left( {{{\left( { - i} \right)}^\gamma } - {e^{i\gamma {\sigma _s}}}} \right)\left( {{{\left( { - i} \right)}^\gamma } - {e^{ - i\gamma {\sigma _s}}}} \right)}  = \]
\[ = {e^{ - 2\mu a + \frac{{\pi k}}{l}a}}{e^{i\frac{{\pi k}}{l}a}}{\left( {\frac{i}{2}} \right)^{m - 1}}{i^q}{M_1}\left( {{{\left( { - i} \right)}^\gamma } - {{\left( { - 1} \right)}^\gamma }} \right)\prod\limits_{s = m + 1}^{2m - 1} {\left( {{{\left( { - 1} \right)}^\gamma } + 1 - 2{{\left( { - i} \right)}^\gamma }\cos \gamma {\sigma _s}} \right)} .\]
Перейдем к вычислению других определителей
\[\det \left( {C_{4m,2m}^ - D_{4m,2m - 1}^ - D_{4m,1}^ - } \right) = \]
\[ = \det \left( {\cos j{\theta _s},\sin j{\theta _s},1, - \cos j{\sigma _t}, - \sin j{\sigma _t},{e^{ - ij\frac{\pi }{2}}}} \right)_{s = \overline {m,2m - 1} ;t = \overline {1,m - 1} }^{j = \overline {0,4m - 1} } = \]
\[ = {\left( {\frac{i}{2}} \right)^{2m - 1}}\det \left( {{e^{ij{\theta _s}}},{e^{ - ij{\theta _s}}},1,{e^{ij{\sigma _k}}},{e^{ - ij{\sigma _k}}},{{\left( { - i} \right)}^j}} \right)_{s = \overline {m,2m - 1} ;k = \overline {1,m - 1} }^{j = \overline {0,4m - 1} } = \]
\[ = {\left( {\frac{i}{2}} \right)^{2m - 1}}{M_2}\left( { - i - 1} \right)\prod\limits_{s = m}^{2m - 1} {\left( { - i - {e^{i{\theta _s}}}} \right)\left( { - i - {e^{ - i{\theta _s}}}} \right)} \prod\limits_{t = 1}^{m - 1} {\left( { - i - {e^{i{\sigma _t}}}} \right)\left( { - i - {e^{ - i{\sigma _t}}}} \right)}  = \]
\[ = {\left( {\frac{i}{2}} \right)^{2m - 1}}{M_2}\left( { - i - 1} \right)\prod\limits_{t = 1}^{m - 1} {2i\cos {\sigma _t}} \prod\limits_{s = m}^{2m - 1} {2i\cos {\theta _s}}  = \]
\[ = {M_2}\left( {i + 1} \right)\prod\limits_{t = 1}^{m - 1} {\cos {\sigma _t}} \prod\limits_{s = m}^{2m - 1} {\cos {\theta _s}}  = {M_3}\left( {i + 1} \right),\]
где
\[{M_2} = \det \left( {{e^{ij{\theta _s}}},{e^{ - ij{\theta _s}}},1,{e^{ij{\sigma _k}}},{e^{ - ij{\sigma _k}}}} \right)_{s = \overline {m,2m - 1} ;k = \overline {1,m - 1} }^{j = \overline {0,4m - 2} } \ne 0.\]
\[{M_3} = {M_2}\prod\limits_{t = 1}^{m - 1} {\cos {\sigma _t}} \prod\limits_{s = m}^{2m - 1} {\cos {\theta _s}}  \ne 0.\]
Аналогично
\[\det \left( {C_{4m,2m}^ - D_{4m,2m - 1}^ - D_{4m,1}^ + } \right) = {M_2}{\left( {\frac{i}{2}} \right)^{2m - 1}}\left( {i - 1} \right)\prod\limits_{t = 1}^{m - 1} {\left( {i - {e^{i{\sigma _k}}}} \right)\left( {i - {e^{ - i{\sigma _k}}}} \right)}  \cdot \]
\[\prod\limits_{s = m}^{2m - 1} {\left( {i - {e^{i{\theta _s}}}} \right)\left( {i - {e^{ - i{\theta _s}}}} \right)}  = {M_2}{\left( {\frac{i}{2}} \right)^{2m - 1}}{\left( { - 2i} \right)^{2m - 1}}\left( {i - 1} \right)\prod\limits_{k = 1}^{m - 1} {\cos {\sigma _k}} \prod\limits_{s = m}^{2m - 1} {\cos {\theta _s}}  = \]
\[ = {M_2}\left( {i - 1} \right)\prod\limits_{k = 1}^{m - 1} {\cos {\sigma _k}} \prod\limits_{s = m}^{2m - 1} {\cos {\theta _s}}  = {M_3}\left( {i - 1} \right).\]
Далее все постоянные не зависящие от $k$ , будем обозначать одной буквой M. Учитывая это имеем
\[{\Delta _{2k}} = M{e^{2a\left( {\alpha  - \mu } \right) + \frac{{\pi k}}{l}a}}{\Delta _{3k}},\]
здесь
\[{\Delta _{3k}} = \left( {{{\left( { - i} \right)}^q}\left( {{i^\gamma } - {{\left( { - 1} \right)}^\gamma }} \right)\prod\limits_{s = m + 1}^{2m - 1} {\left( {{{\left( { - 1} \right)}^\gamma } + 1 - 2{i^\gamma }\cos \gamma {\sigma _s}} \right)} \left( {i + 1} \right){e^{ - \frac{{\pi ka}}{l}i}} - } \right.\]
\[\left. { - {i^q}\left( {{{\left( { - i} \right)}^\gamma } - {{\left( { - 1} \right)}^\gamma }} \right)\prod\limits_{s = m + 1}^{2m - 1} {\left( {{{\left( { - 1} \right)}^\gamma } + 1 - 2{{\left( { - i} \right)}^\gamma }\cos \gamma {\sigma _s}} \right)} \left( {i - 1} \right){e^{i\frac{{\pi k}}{l}a}}} \right).\]
Рассмотрим частные случаи.

1).$\,n = 2m,\,m = 2t + 1,\,\gamma  = 1,\,q = 2p + 1$
\[{\Delta _{3k}} =  - i{\left( { - 1} \right)^p}\prod\limits_{s = m + 1}^{2m - 1} {\left( {2i\cos {\sigma _s}} \right)\left( {\left( {i + 1} \right){{\left( { - 1} \right)}^{m - 1}}\left( {i + 1} \right){e^{ - \frac{{\pi ka}}{l}i}} + \left( {1 - i} \right)\left( {i - 1} \right){e^{i\frac{{\pi k}}{l}a}}} \right)}  = \]
\[ =  - i{\left( { - 1} \right)^p}{\left( { - 4} \right)^t}\prod\limits_{s = m + 1}^{2m - 1} {\cos {\sigma _s}\left( {\left( {i + 1} \right)\left( {i + 1} \right){e^{ - \frac{{\pi ka}}{l}i}} + \left( {1 - i} \right)\left( {i - 1} \right){e^{i\frac{{\pi k}}{l}a}}} \right)}  = \]
\[ = 2{\left( { - 1} \right)^p}{\left( { - 4} \right)^t}\prod\limits_{s = m + 1}^{2m - 1} {\cos {\sigma _s}\sin \left( {\frac{{\pi ka}}{l} + \frac{\pi }{2}} \right)} .\]
Итак
\[{\Delta _{2k}} = M{e^{2a\left( {\alpha  - \mu } \right) + \frac{{\pi k}}{l}a}}\sin \left( {\frac{{\pi ka}}{l} + \frac{\pi }{2}} \right),\]
\[0 \ne M - const,\,\left| M \right| < \infty .\]

2). $n = 2m,\,m = 2t + 1,\,\gamma  = 1,\,q = 2p$
\[{\Delta _{3k}} = {i^{2p}}\prod\limits_{s = m + 1}^{2m - 1} {\left( {2i\cos {\sigma _s}} \right)} \left( {\left( {i + 1} \right)\left( {i + 1} \right){e^{ - \frac{{\pi ka}}{l}i}}} \right.\left. { + \left( {i - 1} \right)\left( {i - 1} \right){e^{i\frac{{\pi k}}{l}a}}} \right) = \]
\[ = 4{i^{2p}}\prod\limits_{s = m + 1}^{2m - 1} {\left( {2i\cos {\sigma _s}} \right)\left( {\sin \frac{{\pi ka}}{l}} \right)} ,\]
итак
\[{\Delta _{2k}} = M{e^{2a\left( {\alpha  - \mu } \right) + \frac{{\pi k}}{l}a}}\sin \frac{{\pi ka}}{l}.\]
Аналогично рассмотрев другие случаи приходим к следующему результату. Основной определитель системы (12) имеет вид
\[{\Delta _{1k}} = M{e^{2a\left( {\alpha  - \mu } \right) + \frac{{\pi k}}{l}a}}\left( {{\Delta _{4k}} + {\Delta _{5k}}} \right),\eqno(13)\]
где
\[0 \ne M - const,\,\,\left| M \right| < \infty ,\,\,\mathop {\lim }\limits_{k \to \infty } {\Delta _{5k}} = 0,\]
\[{\Delta _{4k}} = \left\{ \begin{gathered}
  \sin \left( {\frac{{\pi ka}}{l} + \frac{\pi }{2}} \right),\left\{ {\deg  = 8t + 4,\,\gamma  = 1,\,q = 2p + 1} \right\} \cup \left\{ {\deg  = 8t,\,\gamma  = 1,\,q = 2p} \right\}, \hfill \\
  \sin \frac{{\pi ka}}{l},\left\{ {\deg  = 8t + 4,\,\gamma  = 1,\,q = 2p} \right\} \cup \left\{ {\deg  = 8t,\,\gamma  = 1,\,q = 2p + 1} \right\} \hfill \\
  \sin \left( {\frac{{\pi ka}}{l} + \frac{\pi }{4}} \right),\left\{ {\deg  = 4t,\,\gamma  = 2,q = 2p} \right\}, \hfill \\
  \sin \left( {\frac{{\pi ka}}{l} + \frac{{3\pi }}{4}} \right),\left\{ {\deg  = 4t,\,\gamma  = 2,q = 2p + 1} \right\},\,t = 0,1,..., \hfill \\
\end{gathered}  \right.\eqno(14)\]
deg - порядок уравнения.

В случаи когда $n = 2m + 1,$ получим следующий результат:
\[{\Delta _{1k}} = M{e^{2a\left( {\alpha  - \mu } \right) + \frac{{\pi k}}{l}a}}\left( {{\Delta _{7k}} + {\Delta _{8k}}} \right),\]
где
\[{\Delta _{7k}} = \left\{ \begin{gathered}
  \sin \left( {\frac{{\pi ka}}{l} + \frac{\pi }{4}} \right),\left\{ {\gamma  = 1,q = 2p,deg = 8t+2} \right\}, \hfill \\
  \left\{ {\gamma  = 1,q = 2p + 1,deg = 8t + 6} \right\},\left\{ {\gamma  = 2,q = 2p,deg=4t+2} \right\}; \hfill \\
  \sin \left( {\frac{{\pi ka}}{l} + \frac{{3\pi }}{4}} \right),\left\{ {\gamma  = 1,q = 2p,deg = 8t + 6} \right\}, \hfill \\
  \left\{ {\gamma  = 1,q = 2p + 1,deg = 8t+2} \right\},\left\{ {\gamma  = 2,q = 2p + 1,deg=4t+2} \right\}, \hfill \\
\end{gathered}  \right.\]
\[p,t \in N \cup \left\{ 0 \right\}.\]

\textbf{Теорема 1.} Если существует решение задачи $D$ , то оно единственно только тогда, когда выполнено условие ${\Delta _{1k}} \ne 0$  при всех $k$.

\textbf{Доказательство.} Пусть ${\Delta _{1k}} \ne 0$и граничные условия (4)-(6) однородны, тогда система (11) имеет только тривиальное решение ${u_k}\left( y \right) \equiv 0$ для всех номеров $k$. Значит $u\left( {x,y} \right) = 0$  почти всюду, но т.к. $u(x,y)$ непрерывна в  $\overline \Omega  $ , то и $u\left( {x,y} \right) \equiv 0$  в  $\overline \Omega  $ .

	Пусть теперь при некоторых значениях $a,l,k$ определитель ${\Delta _{1k}} = 0$. Тогда однородная задача (11) будет иметь ненулевое решение, что и будет нарушать единственность задачи D. \textbf{Теорема 1 доказана.}

\begin{center}
\textbf{3.  Существование решения
}\end{center}

Вначале получим некоторые оценки для функций   В дальнейшем, чтобы не увеличивать число обозначений все положительные постоянные не зависящие от к будем обозначать одной буквой $N$. Также будем считать, что $n=2m$ и порядок уравнения равен $4m$. Другой случай рассматривается аналогично.

\textbf{Лемма 1.}  Для функции ${u_k}\left( y \right)$  и её производных, при достаточно больших значениях $k$, справедливы оценки
\[\left| {u_k^{(t)}\left( y \right)} \right| \leqslant N\frac{{{k^t}\sum\limits_{s = 0}^{n - 1} {\left\{ {\left| {{\varphi _{sk}}} \right| + \left| {{\psi _{sk}}} \right|} \right\}} }}{{\left| {{\Delta _{4k}} + {\Delta _{5k}}} \right|}},t = 0,1,...,2n.\]
\textbf{Доказательство.} Нетрудно показать, что
\[\left| {u_k^{\left( t \right)}(y)} \right| \leqslant N{k^t}\left| {{u_k}\left( y \right)} \right|,\]
поэтому достаточно доказать  оценку для $t=0$. Пусть $y>0$, тогда имеем
\[\left| {{u_k}\left( y \right)} \right| \leqslant N\sum\limits_{s = 0}^{2m - 1} {{e^{{\alpha _s}a}}\left( {\left| {c_s^1} \right| + \left| {c_s^2} \right|} \right)} ,\]
\[{e^{{\alpha _0}a}}\left| {c_0^1} \right| = {e^{{\alpha _0}a}}\left| {\frac{{{\Delta _0}}}{{{\Delta _{1k}}}}} \right| \leqslant {e^{{\alpha _0}a}}\frac{{\sum\limits_{s = 0}^{2m - 1} {\left\{ {\left| {{\varphi _{sk}}} \right| + \left| {{\psi _{sk}}} \right|} \right\}}  \cdot O\left( {{e^{2a\left( {\alpha  - \mu } \right) + \frac{{\pi k}}{l}a - {\alpha _0}a}}} \right)}}{{{e^{2a\left( {\alpha  - \mu } \right) + \frac{{\pi k}}{l}a}}\left| {{\Delta _{4k}} + {\Delta _{5k}}} \right|}} \leqslant \]
\[ \leqslant N\frac{{\sum\limits_{s = 0}^{2m - 1} {\left\{ {\left| {{\varphi _{sk}}} \right| + \left| {{\psi _{sk}}} \right|} \right\}} }}{{\left| {{\Delta _{4k}} + {\Delta _{5k}}} \right|}},\]
здесь ${\Delta _0}$ - определитель матрицы, полученный заменой первого столбца основной матрицы системы (12), правой частью системы (12).
Справедливость полученной оценки для других слагаемых показывается аналогично. Случай  $y<0 $рассматривается также. \textbf{Лемма 1 доказана}.

Нужно теперь найти условия, при которых выражение $\left| {{\Delta _{4k}} + {\Delta _{5k}}} \right|$  отделяется от нуля, т.е. начиная с некоторого номера $k \geqslant {k_1} \geqslant 1$  должно выполняться оценка
$$\left| {{\Delta _{4k}} + {\Delta _{5k}}} \right| > \delta  > 0.$$
  Т.к. $\mathop {\lim }\limits_{k \to \infty } {\Delta _{5k}} = 0$  , то достаточно найти условия, при которых   $\left| {{\Delta _{4k}}} \right| > \delta  > 0,\,\forall k.$ Справедлива лемма.

\textbf{Лемма 2.} Для справедливости оценки
$$\left| {{\Delta _{4k}}} \right| \geqslant \delta  > 0,\,\forall k,\eqno(15)$$
достаточно выполнения одного из двух условий: \\
1). ${\Delta _{4k}} = \sin \left( {\frac{{\pi ka}}{l} + \frac{\pi }{4}} \right)$ или ${\Delta _{4k}} = \sin \left( {\frac{{\pi ka}}{l} + \frac{{3\pi }}{4}} \right);$ $\frac{a}{l} \in N,$ либо $\frac{a}{l} = \frac{s}{t}\,\,\,\left( {\frac{a}{l} \notin N} \right),$ $s,t \in N,\,\,\left( {s,t} \right) = 1,$ $t$ не делится на 4\\
2). ${\Delta _{4k}} = \sin \left( {\frac{{\pi ka}}{l} + \frac{\pi }{2}} \right),$ $\frac{a}{l} \in N,$ либо $\frac{a}{l} = \frac{s}{t}\,\,\,\left( {\frac{a}{l} \notin N} \right),$ $s,t \in N,\,\,\left( {s,t} \right) = 1,\,\,\,\left( {t,2} \right) = 1.$

\textbf{Доказательство.} Проверим второе условие (первое проверяется аналогично). Пусть $\frac{a}{l} \in N$, тогда $\left| {\sin \left( {\frac{{\pi ka}}{l} + \frac{\pi }{2}} \right)} \right| = 1$ и все доказано. Пусть теперь $\frac{a}{l} = \frac{s}{t}\,\,\,\left( {\frac{a}{l} \notin N} \right),$ $s,t \in N,\,\,\left( {s,t} \right) = 1,\,\,\,\left( {t,2} \right) = 1,$ тогда $\frac{{ka}}{l} = \frac{{ks}}{t} = {k_1} + \frac{{{k_2}}}{t},$ где ${k_1},{k_2} \in N,\,\,1 \leqslant {k_2} \leqslant t - 1.$ Значит имеем
$$\left| {\sin \left( {\frac{{\pi ka}}{l} + \frac{\pi }{2}} \right)} \right| = \left| {\sin \pi \left( {\frac{{{k_2}}}{t} + \frac{1}{2}} \right)} \right|.$$
Т.к.
$$t \ne 2{k_2} \Rightarrow \frac{{{k_2}}}{t} \ne \frac{1}{2} \Rightarrow \frac{{{k_2}}}{t} + \frac{1}{2} \ne 1 \Rightarrow \left| {\sin \pi \left( {\frac{{{k_2}}}{t} + \frac{1}{2}} \right)} \right| > 0 \Rightarrow $$
\[\left| {\sin \left( {\frac{{\pi ka}}{l} + \frac{\pi }{2}} \right)} \right| \geqslant \delta  = \mathop {\min }\limits_{1 \leqslant {k_2} \leqslant t - 1} \left| {\sin \left( {\frac{{\pi {k_2}}}{t} + \frac{\pi }{2}} \right)} \right| > 0.\]

\textbf{Лемма 2 доказана.}

Отметим, что если $\tau  = \frac{a}{l}$  является иррациональным числом, то не всегда можно отделить знаменатель выражения  от нуля. Но в некоторых случаях можно найти зависимость  «малости» знаменателя  от номера $k$ . Справедлива лемма.

\textbf{Лемма 3.} Если $\tau  = \frac{a}{l} > 0$   является иррациональным алгебраическим числом степени $p \geqslant 2$ , то существует число  $N > 0$ ( не зависимая от $k$) такое , что при всех  $p \in N,\,0 < \varepsilon  < 1$ справедлива оценка
\[\left| {{\Delta _{4k}}} \right| \geqslant \frac{N}{{{k^{1 + \varepsilon }}}}.\eqno(16)\]

\textbf{Доказательство.} Пусть
\[{\Delta _{4k}} = \sin \left( {\pi k\tau  + \frac{\pi }{2}} \right).\]
Используя выпуклость функции $y = \sin x$ на интервале $\left( {0,\frac{\pi }{2}} \right)$ имеем неравенство
\[\left| {\sin x} \right| \geqslant \frac{{2\left| x \right|}}{\pi },\,\left| x \right| \leqslant \frac{\pi }{2},\]
далее имеем
\[\left| {\sin \tau \pi k} \right| = \left| {\sin \pi k\left( {\tau  - \frac{{2m - 1}}{{2k}}} \right)} \right|,\]
где $m\in N$ - произвольно. Теперь подберем $m $ так, чтобы выполнялось неравенство
\[\left| {\tau  - \frac{{2m - 1}}{{2k}}} \right| \leqslant \frac{1}{{2k}},\]
для этого достаточно положить
\[m = \left[ {\tau k} \right] + 1,\]
где $\left[ {\tau k} \right] $- целая часть иррационального числа $\tau k$ . Теперь имеем
\[\left| {\sin \pi k\left( {\tau  - \frac{{2m - 1}}{{2k}}} \right)} \right| \geqslant \frac{2}{\pi }\left| {\pi k\left( {\tau  - \frac{{2m - 1}}{{2k}}} \right)} \right| = 2k\left| {\tau  - \frac{{2m - 1}}{{2k}}} \right|.\]
Известно (Бухштаб А.А. Теория чисел,следствие из теоремы Туэ-Зигель-Рот) что для любого алгебраического числа $\tau $  степени $p \geqslant 2$   и произвольного $0 < \varepsilon  < 1$  найдется $\delta \left( {\tau ,\varepsilon } \right) > 0$  такое, что для любой рациональной дроби $\frac{s}{q}$  выполняется неравенство
\[\left| {\tau  - \frac{s}{q}} \right| \geqslant \frac{{\delta \left( {\tau ,\varepsilon } \right)}}{{{q^{2 + \varepsilon }}}}.\]
Используя это  имеем
\[\left| {\sin \pi k\left( {\tau  - \frac{{2m - 1}}{{2k}}} \right)} \right| \geqslant 2k\left| {\tau  - \frac{{2m - 1}}{{2k}}} \right| \geqslant 2k\frac{\delta }{{{{\left( {2k} \right)}^{2 + \varepsilon }}}} = \frac{N}{{{k^{1 + \varepsilon }}}}.\]
Остальные случаи рассматриваются аналогично. \textbf{Лемма 3 доказана.}

Учитывая вышесказанное, получим условия, при которых ряд
\[u\left( {x,y} \right) = \sum\limits_{k = 1}^\infty  {{u_k}\left( y \right){X_k}\left( x \right)},\eqno(17)\]
является классическим решением поставленной задачи. Формально имеем
\[\left| {D_x^{2n}u\left( {x,y} \right)} \right| \leqslant N\sum\limits_{k = 1}^\infty  {{k^{2n}}} \frac{{\sum\limits_{s = 0}^{n - 1} {\left\{ {\left| {{\varphi _{sk}}} \right| + \left| {{\psi _{sk}}} \right|} \right\}} }}{{\left| {{\Delta _{4k}} + {\Delta _{5k}}} \right|}},\]
теперь, если ${\Delta _{1k}} \ne 0,\,\,\frac{a}{l} -$  удовлетворяет условиям леммы 2 , то имеем оценку
\[\left| {D_x^{2n}u\left( {x,y} \right)} \right| \leqslant N\sum\limits_{k = 1}^\infty  {{k^{2n}}\sum\limits_{s = 0}^{n - 1} {\left\{ {\left| {{\varphi _{ks}}} \right| + \left| {{\psi _{ks}}} \right|} \right\}} },\eqno(18)\]
если ${\Delta _{1k}} \ne 0,\,\,\frac{a}{l} -$  удовлетворяет условиям леммы 3, то
\[\left| {D_x^{2n}u\left( {x,y} \right)} \right| \leqslant N\sum\limits_{k = 1}^\infty  {{k^{2n + 1 + \varepsilon }}\sum\limits_{s = 0}^{n - 1} {\left\{ {\left| {{\varphi _{ks}}} \right| + \left| {{\psi _{ks}}} \right|} \right\}} } ,\,0 < \varepsilon  < 1,\eqno(19)\]
Осталось наложить условия на граничные функции.

\textbf{Теорема 2.} Пусть выполнены следующие условия:\\
1. Определитель ${\Delta _{1k}}$  , системы (12) отличен от нуля;\\
2. Выполнены условия леммы 2;\\
3. ${\varphi _s}\left( x \right),{\psi _s}\left( x \right) \in {C^{2n + 1}}\left[ {0;l} \right],$\\
4. $\varphi _s^{2m}\left( 0 \right) = \varphi _s^{2m}\left( l \right) = \psi _s^{2m}\left( 0 \right) = \psi _s^{2m}\left( l \right) = 0,\,\,s = 0,1,...n - 1,\,\,\,m = 0,1,...n - 1.$\\
Тогда ряд (17) является классическим решением задачи $D.$

\textbf{Доказательство.}
Учитывая (18) имеем
\[\sum\limits_{k = {k_0}}^\infty  {{k^{2n}}\left| {{\varphi _{k0}}} \right|}  = \sum\limits_{k = 1}^\infty  {\frac{1}{k}{k^{2n + 1}}\left| {{\varphi _{k0}}} \right|}  \leqslant \sqrt {\sum\limits_{k = 1}^\infty  {\frac{1}{{{k^2}}}} } \sqrt {\sum\limits_{k = 1}^\infty  {{{\left( {{k^{2n + 1}}{\varphi _{k0}}} \right)}^2}} }  \leqslant N{\left\| {\varphi _0^{(2n + 1)}(x)} \right\|_{{L_{_2}}}}.\]
Сходимость остальных слагаемых доказывается аналогично. \textbf{Теорема 2 доказана.}

\textbf{Теорема 3.} Пусть выполнены следующие условия:\\
1. Определитель ${\Delta _{1k}}$  , системы (12) отличен от нуля;\\
2. Выполнены условия леммы 3;\\
3. ${\varphi _s}\left( x \right),{\psi _s}\left( x \right) \in {C^{2n + 2}}\left[ {0;l} \right],$\\
4. $\varphi _s^{2m}\left( 0 \right) = \varphi _s^{2m}\left( l \right) = \psi _s^{2m}\left( 0 \right) = \psi _s^{2m}\left( l \right) = 0,\,\,s = 0,1,...n - 1,\,\,\,m = 0,1,...n;$\\
5. $\varphi _s^{\left( {2n + 2} \right)}\left( {x + h} \right) - \varphi _s^{\left( {2n + 2} \right)}\left( x \right) = O\left( {{{\left| h \right|}^d}} \right),\,h \to 0,$\\
$\psi _s^{\left( {2n + 2} \right)}\left( {x + h} \right) - \psi _s^{\left( {2n + 2} \right)}\left( x \right) = O\left( {{{\left| h \right|}^d}} \right),\,h \to 0,$\\
$0 < \varepsilon  < d < 1,\,s = 0,1,...,n - 1.$\\
Тогда ряд (17) является классическим решением задачи $D.$\\
\textbf{Доказывается} аналогично теореме 2.

Может оказаться, что выражение  ${\Delta _{1k}} = 0,$ при некоторых значениях $k = {k_1},{k_2},...,{k_p} < {k_0}.$ Тогда для разрешимости задачи (1)-(6), достаточно выполнения условий ${\varphi _{sk}} = \int\limits_0^l {{\varphi _s}\left( x \right){X_k}dx}  = 0,$ ${\psi _{sk}} = \int\limits_0^l {{\psi _s}\left( x \right){X_k}\left( x \right)dx} ,\,s = \overline {0,(n - 1)},$ а само решение будет иметь вид
\[u\left( {x,y} \right) = \sum\limits_{k = 1}^\infty  {\left( {k \ne {k_1},{k_2},...,{k_p}} \right){u_k}\left( y \right){X_k}\left( x \right)}  + \sum\limits_m^{} {\widetilde {{u_m}}\left( y \right){X_k}\left( x \right)} ,\]
где в последней сумме $m$ принимает значения ${k_1},{k_2},...,{k_p},$  функция $\widetilde {{u_m}}\left( y \right)$ - есть ненулевое решение системы (12).

\begin{center}
\textbf{4. Приложение к определению разрешимости краевой задачи для уравнения 4-го порядка
}\end{center}
Рассмотрим уравнение
\[Lu \equiv D_x^4u\left( {x,y} \right) + \left( {\operatorname{sgn} y} \right)\,D_y^4u\left( {x,y} \right) = 0,\eqno(20)\]
в прямоугольной области $\Omega  = \left\{ {\left( {x,y} \right):0 < x < 3, - 1 < y < 1} \right\}.$ Пусть ${\Omega _ + } = \Omega  \cap \left( {y > 0} \right),\,{\Omega _ - } = \Omega  \cap \left( {y < 0} \right).$ Изучим для этого уравнения две краевые задачи.

\textbf{Задача 1.} Найти в области $\Omega $  функцию $u(x,y)$  удовлетворяющую условиям:
\[u \in {C^3}\left( {\overline \Omega  } \right) \cap {C^4}\left( {{\Omega _ + } \cup {\Omega _ - }} \right),\eqno(21)\]
\[Lu\left( {x,y} \right) \equiv 0,\,\left( {x,y} \right) \in {\Omega _ + } \cup {\Omega _ - },\]
\[u\left( {0,y} \right) = u\left( {3,y} \right) = u''\left( {0,y} \right) = u''\left( {3,y} \right) = 0,\,\,\, - 1 \leqslant y \leqslant 1,\eqno(22)\]
\[u\left( {x, - 1} \right) = {\varphi _0}\left( x \right),\,u'\left( {x, - 1} \right) = {\varphi _1}\left( x \right),\,0 \leqslant x \leqslant 3,\eqno(23)\]
\[u\left( {x,1} \right) = {\psi _0}\left( x \right),\,u'\left( {x,1} \right) = {\psi _1}\left( x \right),\,0 \leqslant x \leqslant 3,\eqno(24)\]

\textbf{Задача 2.} Найти в области $\Omega $  функцию $u(x,y)$  удовлетворяющую условиям:
$$(20)-(22),(24),$$
\[u''\left( {x, - 1} \right) = {\varphi _0}\left( x \right),\,u'\left( {x, - 1} \right) = {\varphi _1}\left( x \right),\,0 \leqslant x \leqslant 3.\eqno(25)\]
Относительно ${u_k}\left( y \right)$ имеем обыкновенное диффуравнение
\[u_k^{(4)}\left( y \right) + \operatorname{sgn} y{\left( {\frac{{\pi k}}{3}} \right)^4}{u_k}\left( y \right) = 0,\eqno(26)\]
общее решение (26) имеет вид
\[{u_k}\left( y \right) = {e^{ay}}\left( {{c_1}\cos ay + {c_2}sinay} \right) + {e^{ - ay}}\left( {{c_3}\cos ay + {c_4}sinay} \right),\,(y > 0),\]
здесь
\[a = \frac{{\pi k}}{3}\frac{{\sqrt 2 }}{2},\]
\[{u_k}\left( y \right) = {d_1}{e^{\frac{{\pi k}}{3}y}} + {d_2}\cos \left( {\frac{{\pi k}}{3}y} \right) + {d_3}\sin \left( {\frac{{\pi k}}{3}y} \right) + {d_4}{e^{ - \frac{{\pi k}}{3}y}},\,(y < 0).\]
Для определения неизвестных постоянных получим следующие системы. В случаи задачи 1
\[\left\{ \begin{gathered}
  {c_1}{e^a}\cos a + {c_2}{e^a}sina + {e^{ - a}}\left( {{c_3}\cos a + {c_4}\sin a} \right) = {\psi _0} = A \hfill \\
  {e^a}\left( {{c_1}\cos \left( {a + \frac{\pi }{4}} \right) + {c_2}\sin \left( {a + \frac{\pi }{4}} \right)} \right) +  \hfill \\
   + {e^{ - a}}\left( {{c_1}\cos \left( {a + \frac{{3\pi }}{4}} \right) + {c_2}\sin \left( {a + \frac{{3\pi }}{4}} \right)} \right) = \frac{3}{{\pi k}}{\psi _1} = B \hfill \\
  {d_1}{e^{ - \frac{{\pi k}}{3}}} + {d_2}\cos \left( {\frac{{\pi k}}{3}} \right) - {d_3}\sin \left( {\frac{{\pi k}}{3}} \right) + {d_4}{e^{\frac{{\pi k}}{3}}} = {\varphi _0} = C \hfill \\
  \left( {{d_1}{e^{ - \frac{{\pi k}}{3}}} + {d_2}\sin \left( {\frac{{\pi k}}{3}} \right) + {d_3}\cos \left( {\frac{{\pi k}}{3}} \right) - {d_4}{e^{\frac{{\pi k}}{3}}}} \right) = \frac{3}{{\pi k}}{\varphi _1} = D \hfill \\
  {c_1} + {c_3} - {d_1} - {d_2} - {d_4} = 0, \hfill \\
  \frac{{\sqrt 2 }}{2}{c_1} + \frac{{\sqrt 2 }}{2}{c_2} - \frac{{\sqrt 2 }}{2}{c_3} + \frac{{\sqrt 2 }}{2}{c_4} - {d_1} - {d_3} + {d_4} = 0, \hfill \\
  {c_2} - {c_4} - {d_1} + {d_2} - {d_4} = 0, \hfill \\
   - \frac{{\sqrt 2 }}{2}{c_1} + \frac{{\sqrt 2 }}{2}{c_2} + \frac{{\sqrt 2 }}{2}{c_3} + \frac{{\sqrt 2 }}{2}{c_4} - {d_1} + {d_3} + {d_4} = 0. \hfill \\
\end{gathered}  \right.\eqno(27)\]
В случаи задачи 2
\[\left\{ \begin{gathered}
  {c_1}{e^a}\cos a + {c_2}{e^a}sina + {e^{ - a}}\left( {{c_3}\cos a + {c_4}\sin a} \right) = {\psi _0} = A \hfill \\
  {e^a}\left( {{c_1}\cos \left( {a + \frac{\pi }{4}} \right) + {c_2}\sin \left( {a + \frac{\pi }{4}} \right)} \right) +  \hfill \\
   + {e^{ - a}}\left( {{c_1}\cos \left( {a + \frac{{3\pi }}{4}} \right) + {c_2}\sin \left( {a + \frac{{3\pi }}{4}} \right)} \right) = \frac{3}{{\pi k}}{\psi _1} = B \hfill \\
  {d_1}{e^{ - \frac{{\pi k}}{3}}} - {d_2}\cos \left( {\frac{{\pi k}}{3}} \right) + {d_3}sin\left( {\frac{{\pi k}}{3}} \right) + {d_4}{e^{\frac{{\pi k}}{3}}} = {\left( {\frac{3}{{\pi k}}} \right)^2}{\varphi _0} = C \hfill \\
  \left( {{d_1}{e^{ - \frac{{\pi k}}{3}}} + {d_2}\sin \left( {\frac{{\pi k}}{3}} \right) + {d_3}\cos \left( {\frac{{\pi k}}{3}} \right) - {d_4}{e^{\frac{{\pi k}}{3}}}} \right) = \frac{3}{{\pi k}}{\varphi _1} = D \hfill \\
  {c_1} + {c_3} - {d_1} - {d_2} - {d_4} = 0, \hfill \\
  \frac{{\sqrt 2 }}{2}{c_1} + \frac{{\sqrt 2 }}{2}{c_2} - \frac{{\sqrt 2 }}{2}{c_3} + \frac{{\sqrt 2 }}{2}{c_4} - {d_1} - {d_3} + {d_4} = 0, \hfill \\
  {c_2} - {c_4} - {d_1} + {d_2} - {d_4} = 0, \hfill \\
   - \frac{{\sqrt 2 }}{2}{c_1} + \frac{{\sqrt 2 }}{2}{c_2} + \frac{{\sqrt 2 }}{2}{c_3} + \frac{{\sqrt 2 }}{2}{c_4} - {d_1} + {d_3} + {d_4} = 0. \hfill \\
\end{gathered}  \right.\eqno(28)\]
Учитывая (14), имеем, основные определители систем (27) (${\Delta _1}$) и (28) (${\Delta _2}$) имеют вид
\[{\Delta _1} = M{e^{2a + \frac{{\pi k}}{3}}}\left( {\sin \frac{{\pi k}}{3} + {\Delta _{5k}}} \right),\]
\[{\Delta _2} = M{e^{2a + \frac{{\pi k}}{3}}}\left( {\sin \left( {\frac{{\pi k}}{3} + \frac{\pi }{2}} \right) + {\Delta _{5k}}} \right).\]
Теперь вычислим их с помощью программы Wolfram Mathematica
\[{\Delta _1} =  - 2{e^{ - 2{a_k}}}\left( {2{e^{4{a_k}}} - 4{e^{2{a_k}}}{{\sin }^2}{a_k} + 2} \right) - \]
\[ - 2{e^{ - 2{a_k} - \frac{{\pi k}}{3}}}\left( { - 2{e^{2{a_k}}}\cos \frac{{\pi k}}{3}\left( {\sqrt 2 \left( {{e^{\frac{{2\pi k}}{3}}} - 1} \right)\sin 2{a_k} + {e^{\frac{{2\pi k}}{3}}} + 1} \right)} \right) - \]
 $$- 2{e^{ - 2{a_k} - \frac{{\pi k}}{3}}}\sin \left( {\frac{{\pi k}}{3}} \right)\left( {{e^{4{a_k}}}\left( {\sqrt 2  - 1} \right) - \sqrt 2  - 1} \right)- $$
  $$- 2{e^{ - 2{a_k} + \frac{{\pi k}}{3}}}\sin \left( {\frac{{\pi k}}{3}} \right)\left( {{e^{4{a_k}}} - \sqrt 2  + 1} \right) - $$
\[ - 2\sin \frac{{\pi k}}{3}\left( {2{e^{\frac{{\pi k}}{3}}}{{\sin }^2}{a_k} - 2{{\sin }^2}{a_k} + \sqrt 2 {e^{2{a_k} + \frac{{\pi k}}{3}}}} \right) - \]
\[ - 2{e^{ - 2{a_k} - \frac{{\pi k}}{3}}}\left( { - 2{e^{2{a_k}}}{{\cos }^2}{a_k}\left( {\left( {{e^{\frac{{2\pi k}}{3}}} - 1} \right)\sin \frac{{\pi k}}{3} - 2{e^{\frac{{\pi k}}{3}}}} \right)} \right),\]
видно, что доминирующие члены (самая большая положительная степень экспоненты) совпадают. Далее вычислим по формуле Крамера коэффициент $d_{3}$ (из системы (27))
\[{d_3}\left( k \right) = \frac{1}{{{\Delta _1}}}O\left( {\left( {C + D} \right)\left( {1 + \sqrt 2 } \right){e^{2a + \frac{{\pi k}}{3}}}} \right) \Rightarrow \]
\[\mathop {\lim }\limits_{k = 3l \to \infty } {d_3} = \infty .\]
Значит задача 1 вообще говоря не разрешима методом Фурье.

$${\Delta _2}=2 e^{-2 a-\frac{\pi  k}{3}} \left(-\left(\sqrt{2}-1\right) e^{4 a}+\left(\sqrt{2}-1\right) \left(-e^{\frac{2 \pi  k}{3}}\right)+\sqrt{2}+1\right) \cos \left(\frac{\pi  k}{3}\right)$$
$$+2 e^{-2 a-\frac{\pi  k}{3}} \cos \left(\frac{\pi  k}{3}\right) \left(\left(1+\sqrt{2}\right) e^{4 a+\frac{2 \pi  k}{3}}-2 e^{2 a} \left(e^{\frac{2 \pi  k}{3}}+1\right) \cos (2 a)\right)+$$
$$+2 e^{-2 a-\frac{\pi  k}{3}} \left(4 e^{2 a} \sin \left(\frac{\pi  k}{3}\right) \left(\sqrt{2} \left(e^{\frac{2 \pi  k}{3}}+1\right) \sin (a) \cos (a)+e^{\frac{2 \pi  k}{3}}-1\right)\right),$$
т.к.
\[\cos \frac{{\pi k}}{3} \in \left\{ { - 1, - \frac{1}{2},\frac{1}{2},1} \right\},\]
то
\[{\Delta _2} = O\left( {{e^{2a + \frac{{\pi k}}{3}}}} \right) \Rightarrow \]
задача 2 разрешима единственным образом при условии ${\Delta _2} \ne 0.$

\begin{center}
Литература
\end{center}

1.Бицадзе А.В. Некорректность задачи
Дирихле для уравнений смешанного типа //ДАН СССР. 1953. Т. 122.
№2. С. 167-170.

2. Cannon J. R. A Dirichlet problem for an equation of mixed type
with a discontinuous coefficient //Annali di Matematica Pura ed
Applicata. - 1963. - Т. 61. - №. 1. - С. 371-377.

3. Нахушев А.М. Критерий единственности решения задачи Дирихле для
уравнений смешанного типа в цилиндирической области//
Дифференциальные уравнения. 1970. Т.6. №1. С.190-191.

4. Пташник Б.И. Некорректные граничные задачи для дифференциальных
уравнений с частными производными// Киев, Наукова Думка, 1984, -
С.264.87.

5. Хаджи И.А. Обратная задача для уравнения смешанного типа
с оператором Лаврентьева-Бицадзе // Математические заметки. -
2012.- Т.91.- №.6. - С.908-919.

6. Сабитов К. Б. К теории
начально-граничных задач для уравнения стержней и балок
//Дифференциальные уравнения. - 2017. - Т. 53. - №. 1. - С.
89-100.

 7. Сабитов К. Б. Начальная задача для уравнения колебаний
балки //Дифференциальные уравнения. - 2017. - Т. 53. - №. 5. - С.
665-665.

 8. Сабитов К. Б. Задача Дирихле для уравнений с частными
производными высоких порядков
//Математические заметки. - 2015. - Т. 97. - №. 2. - С. 262-276.

9.Сабитов К.Б. Начально-граничная и обратные задачи для
неоднородного уравнения смешанного параболо-гиперболического
уравнения//Матем.Заметки-2017-Т.102.-№3.-С.415-435.

 10. Нитребич З. М., Пташник Б. Й., Репетило С. М. Задача Діріхле-Неймана для
лінійного гіперболічного рівняння високого порядку зі сталими
коефіцієнтами у смузі //Науковий вісник Ужгородського
університету. Серія: Математика і інформатика. - 2014. - №. 25,№
1. - С. 94-105.

 11. Пташник Б. Й., Репетило С. М. Задача
Діріхле-Неймана у смузі для гіперболічних рівнянь зі сталими
коефіцієнтами //Математичні методи та фізико-механічні поля. -
2013. - №. 56,№ 3. - С. 15-28.

12. Иргашев Б. Ю. Об одной краевой задаче для уравнения высокого
четного порядка //Известия высших учебных заведений. Математика. -
2017. - №. 9. - С. 13-29.

 13. Иргашев Б. Ю. О спектральной задаче
для одного уравнения высокого четного порядка //Известия высших
учебных заведений. Математика. - 2016. - №. 7. - С. 44-54.

 14.Бухштаб А.А. Теория чисел, Просвещение,М.,1966.

\end{document}